# Approximation of Urison operator with operator polynomials of Stancu type

V. L. Makarov[1], I.I. Demkiv[2]

[1] *Instityte of mathematics NAS of Ukraine, Kyiv*

[2] *Lviv polytechnic National University, Lviv, Ukraine*

*October 28, 2011*

Positive polynomial operator that approximates Urison operator, when integration domain is a "regular triangle" is investigated. We obtain Bernstein Polynomials as a particular case.

## 1. INTRODUCTION.

A lot of works are devoted to investigation, generalization and application of Bernstein Polynomials. Let us point to some of them [1] – [3]. Particularly, the work [4] is devoted to application of Bernstein Polynomials for approximation to Urison operator. Here, the linear operator relative to the F is suggested and investigated

$$B_n(F, x(\cdot)) = F(t,0) - \int_0^1 \sum_{k=0}^n \frac{\partial F\left(t, \frac{k}{n} H(\cdot - z)\right)}{\partial z} C_n^k x^k(z)[1-x(z)]^{n-k} dz,$$

that with increasing of its order however precisely approximates Urison operator

$$F(t, x(\cdot)) = \int_0^1 f(t, z, x(z)) dz, \qquad (1)$$

where $H(v)$—is Heaviside function.

In [5] the obtained results have been generalized on the function of two variables, when integration domain is a rectangle.

The purpose of this work is investigation of linear operator polynomial that approximates Urison operator, when integration domain is a "regular triangle". As the basis of this polynomial D.D. Stancu polynomials are suggested (see [2], [3]), and serve for obtaining Bernstein Polynomials as a particular case.

## 2. THE CASE OF ONE VARIABLE

Let us consider Urison operator (1) with unknown kernel $f(t, z, x(z))$, its properties can be judged only by its influence on any functions $x(z)$ from certain class. In technique such situation is sometimes called "grey box".

The problem lies in the development of simple polynomial approximation of operator $F$, that with the increasing of its order could approximate $F$ with high accuracy.

Let the continuous interpolation conditions are set

$$F(t, x_i(\cdot, \xi)) = \int_0^1 f(t, z, x_i(z, \xi)) dz, \qquad i = \overline{0, n},$$

where

$$x_i(z,\xi) = \frac{i}{n} H(z-\xi), \qquad \xi \in [0,1], \quad i = \overline{0,n}.$$

It is necessary to develop the approximation to Urison operator (1) with unknown function $f(t,z,x(z))$ by means of noted function $F(t,x_i(\cdot,\xi\cdot))$. As the basis of such approximation we will take a linear positive operator polynomial, investigated in [2], which approximates the function of one variable, defined and bounded in the interval $[0,1]$. We have

$$P_n^{[\alpha]}(F,x(\cdot)) = \int_0^1 \sum_{k=0}^n C_n^k v_n^{(k)}(x(z),\alpha) f\left(t,z,\frac{k}{n}\right) dz = \int_0^1 p_n^{[\alpha]}(f(t,z,\cdot);x(z)) dz, \qquad (2)$$

where $v_n^{(k)}(x(z),\alpha) = \dfrac{\prod_{k_1=0}^{k-1}(x(z)+k_1\alpha)\prod_{k_2=0}^{n-k-1}(1-x(z)+k_2\alpha)}{\prod_{k_3=0}^{n-1}(1+k_3\alpha)}$, $\alpha$ is nonnegative parameter that can depend only on n.

Note, that the node $x_0(z,\xi) = 0$ and continuous node $x_n(z,\xi) = H(z-\xi)$ will be the interpolation nodes. If $\alpha = -\dfrac{1}{n}$, then all functions , $i = 0,1,...,n$, $\xi \in [0,1]$ will be continuous interpolation nodes. Really, in this case formula (2) acquires the form

$$P_n^{\left[-\frac{1}{n}\right]}(F,x(\cdot)) = \int_0^1 \sum_{k=0}^n C_n^k v_n^{(k)}\left(x(z),-\frac{1}{n}\right) f\left(t,z,\frac{k}{n}\right) dz = \int_0^1 p_n^{\left[-\frac{1}{n}\right]}(f(t,z,\cdot);x(z)) dz.$$

Let substitute the function $x_i(z,\xi) = \dfrac{i}{n} H(z-\xi)$ into it, then we obtain

$$P_n^{\left[-\frac{1}{n}\right]}\left(F,\frac{i}{n}H(\cdot-\xi)\right) = \int_0^\xi \sum_{k=0}^n C_n^k v_n^{(k)}\left(0,-\frac{1}{n}\right) f\left(t,z,\frac{k}{n}\right) dz + \int_\xi^1 \sum_{k=0}^n C_n^k v_n^{(k)}\left(\frac{i}{n},-\frac{1}{n}\right) f\left(t,z,\frac{k}{n}\right) dz =$$

$$= \int_0^\xi f(t,z,0) dz + \int_\xi^1 f\left(t,z,\frac{i}{n}\right) dz = F(t,\frac{i}{n}H(\cdot-\xi)), \quad i = 0,1,...,n, \quad \xi \in [0,1].$$

Here we used the equalities

$$v_n^{(k)}\left(\frac{i}{n},-\frac{1}{n}\right) = \delta_{i,k}, \quad i,k = 0,1,...,n,$$

where $\delta_{i,k}$ – the symbol of Kronecker. Because of uniqueness the interpolation polynomial $p_n^{[-\frac{1}{n}]}(f(t,z,\cdot);x(z))$ coincides with interpolaion Lagrange polynomial.

The following identities have been proved in [2]

$$p_n^{[\alpha]}(1;x) = 1, \quad p_n^{[\alpha]}(t;x) = t, \quad p_n^{[\alpha]}(t^2;x) = \frac{1}{1+\alpha}\left[\frac{x(1-x)}{n} + x(x+\alpha)\right] \qquad (3)$$

In (2) the unknown functions $f\left(t, z, \dfrac{k}{n}\right)$, $i = \overline{0, n}$ are involved. To define them we use the results of work [6]. Then we will have

$$f\left(t, z, \frac{k}{n}\right) = -\frac{\partial F\left(t, \dfrac{k}{n} H(\cdot - z)\right)}{\partial z} + f(t, z, 0). \qquad (4)$$

Let us use (4) and reduce (2) to the form

$$P_n^{[\alpha]}(F, x(\cdot)) = F(t, 0) - \int_0^1 \sum_{k=0}^n C_n^k v_n^k(x(z), \alpha) \frac{\partial F\left(t, \dfrac{k}{n} H(\cdot - z)\right)}{\partial z} dz,$$

that is used in the concrete calculation, but while proving the convergence theorem we will take into the account formula (2). Let's suppose that the condition holds

$$f(t) \in C([0,1]^3). \qquad (5)$$

Further we will need the theorem of P.P. Korovkin [7].
Let

$$L_n(f; x) = \int_a^b f(t) d\phi_n(x, t), \quad n = 1, 2, \ldots$$

is the sequence of linear operators, defined for $f(t) \in C[a, b]$, where $\phi_n(x, t)$ for every $n$ and for every fixed $x$ is a function of bounded variation of variable $t$ in [a,b]. Then the following statement is valid:

*If in $[a, b]$ $L_n(t^i; x)$ converges uniformly to $t^i$, $i = 0, 1, 2$, then $L_n(f; x)$ converges uniformly to $f(x)$ for any $f(t) \in C[a, b]$.*

The following theorem is true.

**Theorem 1.** *Let the Urison operator (1) is such, that the function $f(t, z, x)$ satisfies the condition (5), and let the operator (1) is considered in the compact $\Phi \subset C[0,1]$ and $0 \leq \alpha = \alpha(n) \to 0$ as $n \to \infty$, then the sequence of operators $\{P_n^\alpha(F, x(\cdot))\}$ converges to $F(t, x(\cdot))$ uniformly with respect to $x(t) \in \Phi$, where $\Phi = \{x(z) \in C[0,1] : 0 \leq x(t) \leq 1\}$.*

**Proof.**
For every fixed $t, z \in [0,1]$ operator polynomial $P_n^{[\alpha]}(f(t, z, \cdot); x(z))$ from (2), with respect to Korovkin P.P. theorem converges to $f(t, z, x(z))$ anywhere on the compact $\Phi \subset C[0,1]$, if $0 \leq \alpha = \alpha(m) \to 0$ as $m \to \infty$. Then, from the evident generalization of theorem 1, p. 506 from [8] (Boundary conversion under integral sign. Chapter. Integrals, that depend on the parameter), in our case the theorem 1 statement follows.

### 3. ESTIMATE OF THE APPROXIMATION ORDER IN CASE OF ONE VARIABLE
Next we will use the modulus of continuity, defined as follows

$$\omega(f, \delta) = \omega(\delta) = \max_{t, z \in [0,1]} \sup_{|x - \tilde{x}| \leq \delta} |f(t, z, x) - f(t, z, \tilde{x})|,$$

where $\delta$ does not depend on $x$, $\tilde{x}$.

**Theorem 2.** *Let Urison operator (1) is such, that the function $f(t, z, x)$ satisfies condition (5), and let operator (1) is considered on compact $\Phi \subset C[0,1]$ and $\alpha \geq 0$, then*

$$\left| F(t, x(\cdot)) - P_n^{[\alpha]}(F, x(\cdot)) \right| \leq \frac{3}{2} \omega\left( \sqrt{\frac{1+\alpha n}{n+\alpha n}} \right). \tag{6}$$

**Proof.** Whereas, on the compact $\Phi \subset C[0,1]$ we have $C_n^k v_n^{(k)}(x(z), \alpha) \geq 0$ and $P_m^{[\alpha]}(1, x(\cdot)) = 1$, then it might be written

$$\left| F(t, x(\cdot)) - P_n^{[\alpha]}(F, x(\cdot)) \right| \leq \int_0^1 \sum_{k=0}^n \left| f(t, z, x(z)) - f\left(t, z, \frac{k}{n}\right) \right| C_n^k v_n^{(k)}(x(z), \alpha) dz.$$

Let us use the following module continuity properties

$$|g(x) - g(\tilde{x})| \leq \omega(|x - \tilde{x}|), \quad \omega(\lambda \delta) \leq (1 + \lambda) \omega(\delta), \ \lambda > 0. \tag{7}$$

Then we obtain

$$\left| f(t, z, x(z)) - f\left(t, z, \frac{k}{n}\right) \right| \leq \omega\left( f(t, z, \cdot); \left| x(z) - \frac{k}{n} \right| \right) \leq (1 + \frac{1}{\delta} \left| x(z) - \frac{k}{n} \right|) \omega(f(t, z, \cdot); \delta) =$$

$$= (1 + \frac{1}{\delta} \left| x(z) - \frac{k}{n} \right|) \omega(\delta).$$

Hence

$$\left| F(t, x(\cdot)) - P_n^{[\alpha]}(F, x(\cdot)) \right| \leq \int_0^1 \sum_{k=0}^n \left( 1 + \frac{1}{\delta} \left| x(z) - \frac{k}{n} \right| \right) \omega(\delta) C_n^k v_n^{(k)}(x(z), \alpha) dz =$$

$$= \left( 1 + \frac{1}{\delta} \int_0^1 \sum_{k=0}^n C_n^k v_n^{(k)}(x(z), \alpha) \left| x(z) - \frac{k}{n} \right| dz \right) \omega(\delta). \tag{8}$$

Now taking into account both Cauchy – Schwarz inequality and identities (3) we write

$$\int_0^1 \sum_{k=0}^n C_n^k v_n^{(k)}(x(z), \alpha) \left| x(z) - \frac{k}{n} \right| dz \leq \int_0^1 \left[ \sum_{k=0}^n C_n^k v_n^{(k)}(x(z), \alpha) \left| x(z) - \frac{k}{n} \right|^2 \right]^{\frac{1}{2}} dz =$$

$$= \int_0^1 \left[ x^2(z) - 2x(z) \sum_{k=0}^n C_n^k v_n^{(k)}(x(z), \alpha) \frac{k}{n} + \sum_{k=0}^n C_n^k v_n^{(k)}(x(z), \alpha) \left( \frac{k}{n} \right)^2 \right]^{\frac{1}{2}} dz \leq$$

$$\leq \max_{x(z) \in [0,1]} \int_0^1 \left[ \frac{x(z)(1 - x(z))}{n} \frac{1 + \alpha n}{1 + \alpha} \right]^{\frac{1}{2}} dz \leq \frac{1}{2} \sqrt{\frac{1 + \alpha n}{n + \alpha n}}.$$

Using this, we write (8) in the form

$$\left| F(t, x(\cdot)) - P_n^{[\alpha]}(F, x(\cdot)) \right| \leq \left( 1 + \frac{1}{2\delta} \sqrt{\frac{1 + \alpha n}{n + \alpha n}} \right) \omega(\delta).$$

To obtain the statement of Theorem 1, we will choose $\delta = \sqrt{\dfrac{1+\alpha n}{n+\alpha n}}$. So we have proved the theorem.

In the case when $\alpha = 0$, (6) transforms into the following inequality for the operator of Bernstein type

$$|F(t,x(\cdot)) - B_n(F,x(\cdot))| \le \frac{3}{2}\sqrt{\frac{1}{n}}.$$

## 4. ASYMPTOTIC ESTIMATE OF APPROXIMATION ERROR IN CASE OF ONE VARIABLE

Let us determine an asymptotic estimate for the error

$$R_n^{[\alpha]}(F, x(\cdot)) = F(t, x(\cdot)) - P_n^{[\alpha]}(F, x(\cdot)).$$

The following theorem is true.

**Theorem 3.** *Let Urison operator (1) is such, that the function $f(t,z,x(\cdot))$ satisfies condition (5) and has the second continuous derivative with respect to $x \in [0,1]$. Then we have asymptotic formula for the remainder*

$$R_n^{[\alpha]}(F, x(\cdot)) = -\frac{1}{2} \cdot \frac{1+\alpha n}{1+\alpha} \int_0^1 \frac{x(z)(1-x(z))}{n} \frac{\partial^2}{\partial x^2} f(t,z,x(z))dz + \int_0^1 \varepsilon_n^{[\alpha]}(x(z))dz, \qquad (9)$$

where $n\varepsilon_n^{[\alpha]}(x(z))$ tends to zero as $n \to \infty$.

*Proof.* From (1) and (2) taking into account (3) one can write

$$R_n^{[\alpha]}(F, x(\cdot)) = \int_0^1 \sum_{k=0}^n C_n^k v_n^k(x(z),\alpha)\left[f(t,z,x(z)) - f(t,z,\frac{k}{n})\right]dz.$$

Let us substitute Taylor series at the point $x(z)$ with the remainder in the integral form

$$f\left(t,z,\frac{k}{n}\right) = f(t,z,x(z)) + \left(\frac{k}{n} - x(z)\right)\frac{\partial}{\partial x}f(t,z,x(z)) + \int_{x(z)}^{k/n}\left(\frac{k}{n} - s\right)\frac{\partial^2}{\partial s^2}f(t,z,s)ds$$

for $f\left(t,z,\dfrac{k}{n}\right)$.

Let us add and subtract $\dfrac{\partial^2}{\partial x^2}f(t,z,x(z))$ in the integral error. We obtain

$$R_n^{[\alpha]}(F, x(\cdot)) = -\int_0^1 \sum_{k=0}^n C_n^k v_n^k(x(z),\alpha)\left(\frac{k}{n} - x(z)\right)\frac{\partial}{\partial x}f(t,z,x(z))dz -$$

$$-\int_0^1 \sum_{k=0}^n C_n^k v_n^k(x(z),\alpha)\left[\int_{x(z)}^{k/n}\left(\frac{k}{n} - s\right)\left[\frac{\partial^2}{\partial s^2}f(t,z,s) - \frac{\partial^2}{\partial x^2}f(t,z,x(z))\right]ds\right]dz -$$

$$-\int_0^1 \sum_{k=0}^n C_n^k v_n^k(x(z),\alpha)\left[\int_{x(z)}^{k/n}\left(\frac{k}{n}-s\right)\frac{\partial^2}{\partial x^2}f(t,z,x(z))ds\right]dz.$$

Let us use the identities (3). We will get

$$R_n^{[\alpha]}(F,x(\cdot)) = \int_0^1\left[-\frac{1}{2}\cdot\frac{1+\alpha n}{1+\alpha}\cdot\frac{x(z)(1-x(z))}{n}\frac{\partial^2}{\partial x^2}f(t,z,x(z))-\right.$$

$$\left.-\sum_{k=0}^n C_n^k v_n^{(k)}(x(z),\alpha)\int_{x(z)}^{\frac{k}{n}}\left(\frac{\partial^2}{\partial s^2}f(t,z,s)-\frac{\partial^2}{\partial x^2}f(t,z,x(z))\right)\left(\frac{k}{n}-s\right)ds\right]dz =$$

$$= -\frac{1}{2}\cdot\frac{1+\alpha n}{1+\alpha}\int_0^1\frac{x(z)(1-x(z))}{n}\frac{\partial^2}{\partial x^2}f(t,z,x(z))dz + \int_0^1 \varepsilon_n^{[\alpha]}(x(z))dz.$$

Now, let us estimate $\varepsilon_n^{[\alpha]}(x(z))$.

$$\varepsilon_n^{[\alpha]}(x(z)) \leq \sum_{k=0}^n C_n^k v_n^{(k)}(x(z),\alpha)\left|\int_{x(z)}^{\frac{k}{n}}\omega(\frac{\partial^2}{\partial x^2}f(t,z,x);|s-x(z)|)\left|\frac{k}{n}-s\right|ds\right| \leq$$

$$\leq \sum_{k=0}^n C_n^k v_n^{(k)}(x(z),\alpha)\left|\int_{x(z)}^{\frac{k}{n}}\left(1+\frac{1}{\delta}|s-x(z)|\right)\omega(\frac{\partial^2}{\partial x^2}f(t,z,x);\delta)\left|\frac{k}{n}-s\right|ds\right| \leq$$

$$\leq \sum_{k=0}^n C_n^k v_n^{(k)}(x(z),\alpha)\left|\int_{x(z)}^{\frac{k}{n}}\left(1+\frac{1}{\delta}|s-x(z)|\right)\left|\frac{k}{n}-s\right|ds\right|\omega(\frac{\partial^2}{\partial x^2}f;\delta) \leq$$

$$\leq \sum_{k=0}^n C_n^k v_n^{(k)}(x(z),\alpha)\left(\frac{1}{2}\left(\frac{k}{n}-x(z)\right)^2+\frac{1}{6\delta}\left|\frac{k}{n}-x(z)\right|^3\right)\omega(\frac{\partial^2}{\partial x^2}f;\delta) \leq$$

$$\leq \left\{\sum_{k=0}^n C_n^k v_n^{(k)}(x(z),\alpha)\left(\frac{k}{n}-x(z)\right)^4\right\}^{1/2}\left\{\sum_{k=0}^n C_n^k v_n^{(k)}(x(z),\alpha)\left(\frac{1}{2}+\frac{1}{6\delta}\left|\frac{k}{n}-x(z)\right|\right)^2\right\}^{1/2}\omega(\frac{\partial^2}{\partial x^2}f;\delta) \leq$$

$$\leq \left\{\sum_{k=0}^n C_n^k v_n^{(k)}(x(z),\alpha)\left(\frac{k}{n}-x(z)\right)^4\right\}^{1/2}\left\{\sum_{k=0}^n C_n^k v_n^{(k)}(x(z),\alpha)\left(\frac{1}{2}+\frac{1}{18\delta^2}\left|\frac{k}{n}-x(z)\right|^2\right)\right\}^{1/2}\omega(\frac{\partial^2}{\partial x^2}f;\delta).$$

Whereas, the functions $v_n^{(k)}(x(z),\alpha)$ is continuous with respect to $\alpha$, we consider the limit of the last inequality, as $\alpha \to 0$. So, we have

$$\varepsilon_n^{[0]}(x(z)) \leq \left\{\sum_{k=0}^n C_n^k x(z)^k(1-x(z))^{n-k}\left(\frac{k}{n}-x(z)\right)^4\right\}^{1/2} \times$$

$$\times \left\{\sum_{k=0}^n C_n^k x(z)^k(1-x(z))^{n-k}\left(\frac{1}{2}+\frac{1}{18\delta^2}\left|\frac{k}{n}-x(z)\right|^2\right)\right\}^{1/2}\omega(\frac{\partial^2}{\partial x^2}f;\delta) \leq$$

$$\leq \sqrt{\frac{1}{n^3}\left(\frac{3n}{16}-\frac{1}{8}\right)\left(\frac{1}{2}+\frac{1}{72\delta^2 n}\right)}^{1/2} \omega(\frac{\partial^2}{\partial x^2}f;\delta) \leq \frac{1}{4n}\sqrt{\frac{37}{24}}\omega(\frac{\partial^2}{\partial x^2}f;\frac{1}{\sqrt{n}}) \qquad (10)$$

To obtain this inequality has been done the following
1) in the first brace we have used the identity

$$\sum_{k=0}^{n} C_n^k x(z)^k (1-x)^{n-k} \left(\frac{k}{n}-x\right)^4 = \frac{x(x-1)}{n^3}\left(3xn(x-1)-6x^2+6x-1\right) \qquad (11)$$

and the fact, that its right side gains maximum $\frac{1}{n^3}\left(\frac{3n}{16}-\frac{1}{8}\right)$ with respect to $x$ over interval [0,1], at $x=\frac{1}{2}$.

2) in the second brace we have used the identities (3).

3) we have chosen $\delta = \frac{1}{\sqrt{n}}$.

Let us denote $n\varepsilon_n^{[\alpha]}(x(z)) = \frac{1}{4}\sqrt{\frac{37}{24}}\omega(\frac{\partial^2}{\partial x^2}f;\frac{1}{\sqrt{n}})$.

From (10) follows, that if $n$ is big enough the theorem will be valid.

The theorem is proved.

In a particular case, when $\alpha = 0$, (9) turns into the asymptotic formula for Bernstein polynomial

$$B_n(F,x(\cdot)) = F(t,x(\cdot)) + \frac{1}{2}\cdot\int_0^1 \frac{x(z)(1-x(z))}{n}\frac{\partial^2}{\partial x^2}f(t,z,x(z))dz - \int_0^1 \frac{\varepsilon_n^{[0]}(x(z))}{n}dz.$$

If $\alpha = -\frac{1}{n}$, then, as it has been noticed previously, the polynomial $p_n^{[\alpha]}(f(t,z,\cdot);x(z))$ coincides with Lagrange interpolation polynomial. Thus, under condition $\frac{\partial^{n+1}}{\partial x^{n+1}}f(t,z,x) \in C([0,1]^3)$, the following formula is valid

$$P_n^{\left[-\frac{1}{n}\right]}(F,x(\cdot)) = F(t,x(\cdot)) + \frac{1}{(n+1)!}\int_0^1 \prod_{k=0}^{n}(x(z)-\frac{k}{n})\frac{\partial^{n+1}}{\partial x^{n+1}}f(t,z,\xi(x(z)))dz, \qquad (12)$$

where $x(z)) \in (0,1)$ is some intermediate point.

*Example 1.* Let's consider Urison operator

$$F(t,x(\cdot)) = \int_0^1 \sin(z^2 x(z))dz, \quad x(z) \in \Phi.$$

and construct operator polynomial (2) for it

$$P_n^{[\alpha]}(F, x(\cdot)) = \int_0^1 \sum_{k=0}^n C_n^k \frac{\prod_{k_1=0}^{k-1}(x(z)+k_1\alpha)\prod_{k_2=0}^{n-k-1}(1-x(z)+k_2\alpha)}{\prod_{k_3=0}^{n-1}(1+k_3\alpha)} \sin\left(z^2 \frac{k}{n}\right) dz.$$

For instance, let's choose, $\alpha = \frac{1}{n}$, $\alpha = -\frac{1}{n}$ and $x(z) = 1 - z^2/(1+z^2)$. For calculations we'll use Maple. Let's set **Digits:=256.** The results we write into table 1, where
$\Delta_1 = |F(t, x(\cdot)) - P_n^{[1/n]}(F, x(\cdot))|$, $\Delta_2 = |F(t, x(\cdot)) - B_n(F, x(\cdot))|$, $\Delta_3 = |F(t, x(\cdot)) - P_n^{[-1/n]}(F, x(\cdot))|$.

Table 1

| n | $\Delta_1$ | $C_1 = n*\Delta_1$ | $\Delta_2$ | $C_2 = n*\Delta_2$ | $\Delta_3$ |
|---|---|---|---|---|---|
| 4 | 0.00359457 | 0.0143783 | 0.00229443 | 0.0091777 | 3.9117323e-7 |
| 8 | 0.00203442 | 0.016275 | 0.0011621 | 0.009296 | 8.4629711e-14 |
| 16 | 0.0010908 | 0.0174528 | 0.00058493 | 0.0093589 | 6.0363698e-29 |
| 32 | 0.000566202 | 0.0181185 | 0.000293464 | 0.0093909 | 2.1917865e-63 |
| 64 | 0.00028866 | 0.0184741 | 0.000146985 | 0.009407 | 3.4561993e-141 |

We see, that the following inequalities hold

$$\left|F(t, x(\cdot)) - P_n^{\left[\frac{1}{n}\right]}(F, x(\cdot))\right| \leq \frac{0{,}019}{n}, \qquad |F(t, x(\cdot)) - B_n(F, x(\cdot))| \leq \frac{0{,}095}{n}.$$

Let's use formula (12). Then, using the reasoning from [9], p.95 it is easy to convince that the following estimate is valid

.

Its right side tends to zero very fast. Thus, if $n = 4$ it is already less than .

## 5. THE CASE OF TWO VARIABLES

Consider Urison operator in case of two variables

$$F(t, x(\cdot), y(\cdot)) = \iint_{\Delta_2} f(t, z_1, z_2, x(z_1), y(z_2)) dz_1 dz \qquad (13)$$

with unknown kernel $f(t, z_1, z_2, x(z_1), y(z_2))$ and integration domain

$$\Delta_2 = \{(z_1, z_2): z_1 \geq 0; z_2 \geq 0; z_1 + z_2 \leq 1\}.$$

Besides, let the following inequalities hold $x(z_1) \geq 0$, $y(z_2) \geq 0$, $x(z_1) + y(z_2) \leq 1$.
Let the interpolation conditions are set

$$F(t, x_i(\cdot), y_j(\cdot)) = \iint_{\Delta_2} f(t, z_1, z_2, x_i(z_1), y_j(z_2)) dz_1 dz,$$

Where

$$x_i(z_1, \xi_1) = \frac{i}{m} H(z_1 - \xi_1), \qquad \xi_1 \in [0,1], \qquad i = \overline{0, m},$$

$$y_j(z_2, \xi_2) = \frac{j}{m} H(z_2 - \xi_2), \quad \xi_2 \in [0,1], \qquad j = \overline{0, m},$$

$H(\gamma)$ – is a Heaviside function, $\xi_1 + \xi_2 \leq 1$, $\quad 0 \leq x_i(z_1, \xi_1) + y_j(z_2, \xi_2) \leq 1$, $\quad i + j \leq m$.

The problem lies in the development of operator polynomial approximation to operator (13) by means of defined function $F(t, x_i(\cdot), y_j(\cdot))$ which can approximate $F(t, x(\cdot), y(\cdot))$ as precisely as possible with the increasing of its order. As the bases we'll take a linear positive polynomial operator from [3], which approximates the two variable function $f(x, y)$, that is bounded in the regular triangle

$$\Delta = \{(x, y): x \geq 0; y \geq 0; x + y \leq 1\}.$$

We gain the next operator polynomial

$$P_m^{[\alpha]}(F, x(\cdot), y(\cdot)) = \iint\limits_{\Delta_2} \sum_{0 \leq i+j \leq m} w_m^{(i,j)}(x(z_1), y(z_2), \alpha) f\left(t, z_1, z_2, \frac{i}{m}, \frac{j}{m}\right) dz_1 dz_2 =$$
$$= \iint\limits_{\Delta_2} p_m^{[\alpha]}(t, z_1, z_2, x(\cdot), y(\cdot)) dz_1 dz_2, \qquad (14)$$

where $w_m^{(i,j)}(x(z_1), y(z_2), \alpha) = C_{i,j}(m) v_m^{(i,j)}(x(z_1), y(z_2), \alpha)$, $C_{i,j}(m) = \dfrac{m!}{i! \, j! \, (m-i-j)!}$,

$$v_m^{(i,j)}(x(z_1), y(z_2), \alpha) = \frac{\prod\limits_{k_1=0}^{i-1}(x(z_1) + k_1 \alpha) \prod\limits_{k_2=0}^{j-1}(y(z_2) + k_2 \alpha) \prod\limits_{k_3=0}^{m-i-j-1}(1 - x(z_1) - y(z_2) + k_3 \alpha)}{\prod\limits_{k_4=0}^{m-1}(1 + k_4 \alpha)},$$

$\alpha$ is a nonnegative parameter that depends only on m.

Note, that from $w_m^{(i,j)}(x(z_1), y(z_2), \alpha) \geq 0$ y $\Delta_2$ follows, that linear operator (14) is positive in $\Delta_2$ in the sense, that if $f\left(t, z_1, z_2, \dfrac{i}{m}, \dfrac{j}{m}\right) \geq 0$ in $\Delta_2$, then $P_m^{[\alpha]}(F, x(\cdot), y(\cdot)) \geq 0$ in $\Delta_2$.

In the case $\alpha = 0$ operator (14) transforms into Bernstein operator in the form

$$B_m(F, x(\cdot), y(\cdot)) = \iint\limits_{\Delta_2} \sum_{0 \leq i+j \leq m} \frac{m!}{i! \, j! \, (m-i-j)!} x^i(z_1) y^j(z_2) (1 - x(z_1) - y(z_2))^{m-i-j} f\left(t, z_1, z_2, \frac{i}{m}, \frac{j}{m}\right) dz_1 dz_2.$$

Lets show, similarly to one dimensional case, that if $\alpha = -\dfrac{1}{m}$ the nodes

$$\left\{\frac{i}{m} H(z_1 - \xi_1), \frac{j}{m} H(z_2 - \xi_2)\right\}, \; 0 \leq i + j \leq m, \; \{z_1, z_2\} \in \Delta_2,$$

$$\left\{\frac{i}{m} H(z_1 - \xi_1), 0\right\}, \left\{0, \frac{j}{m} H(z_2 - \xi_2)\right\}, \; z_1, z_2, \xi_1, \xi_2 \in [0,1]$$

will be continuous interpolation nodes. For the last two nodes this follows from one-dimensional case. Then we have

$$P_m^{\left[-\frac{1}{m}\right]}\left(F, \frac{k}{m}H(\cdot-\xi_1), \frac{l}{m}H(\cdot-\xi_2)\right) =$$

$$\iint_{\Delta_2} \sum_{0\le i+j\le m} w_m^{(i,j)}\left(\frac{k}{m}H(z_1-\xi_1), \frac{l}{m}H(z_2-\xi_2), -\frac{1}{m}\right) f\left(t, z_1, z_2, \frac{i}{m}, \frac{j}{m}\right) dz_1 dz_2 =$$

$$= \int_0^{\xi_1} \int_0^{1-\xi_1} \sum_{0\le i+j\le m} w_m^{(i,j)}\left(0, 0, -\frac{1}{m}\right) f\left(t, z_1, z_2, \frac{i}{m}, \frac{j}{m}\right) dz_2 dz_1 +$$

$$+ \int_0^{\xi_1} \int_{1-\xi_1}^{1-z_1} \sum_{0\le i+j\le m} w_m^{(i,j)}\left(0, \frac{l}{m}, -\frac{1}{m}\right) f\left(t, z_1, z_2, \frac{i}{m}, \frac{j}{m}\right) dz_2 dz_1 +$$

$$+ \int_{\xi_1}^{1-\xi_2} \int_0^{1-z_1} \sum_{0\le i+j\le m} w_m^{(i,j)}\left(\frac{k}{m}, \frac{l}{m}, -\frac{1}{m}\right) f\left(t, z_1, z_2, \frac{i}{m}, \frac{j}{m}\right) dz_2 dz_1 +$$

$$+ \int_{1-\xi_2}^{1} \int_0^{1-z_1} \sum_{0\le i+j\le m} w_m^{(i,j)}\left(\frac{k}{m}, 0, -\frac{1}{m}\right) f\left(t, z_1, z_2, \frac{i}{m}, \frac{j}{m}\right) dz_2 dz_1.$$

Or, taking into account relation [3] $w_m^{(i,j)}\left(\frac{k}{m}, \frac{l}{m}, -\frac{1}{m}\right) = \delta_{i,k}\,\delta_{j,l}$, $i,j,k,l = 0,1,\ldots,m$, finally we obtain

$$P_m^{\left[-\frac{1}{m}\right]}\left(F, \frac{k}{m}H(\cdot-\xi_1), \frac{l}{m}H(\cdot-\xi_2)\right) = \int_0^{\xi_1}\int_0^{1-\xi_1} f(t, z_1, z_2, 0, 0) dz_2 dz_1 + \int_0^{\xi_1}\int_{1-\xi_1}^{1-z_1} f\left(t, z_1, z_2, 0, \frac{l}{m}\right) dz_2 dz_1 +$$

$$+ \int_{\xi_1}^{1-\xi_2}\int_0^{1-z_1} f\left(t, z_1, z_2, \frac{k}{m}, \frac{l}{m}\right) dz_2 dz_1 + \int_{1-\xi_2}^{1}\int_0^{1-z_1} f\left(t, z_1, z_2, \frac{k}{m}, 0\right) dz_2 dz_1 = F(t, \frac{k}{m}H(\cdot-\xi_1), \frac{l}{m}H(\cdot-\xi_2)).$$

In (14) unknown functions $f\left(t, z_1, z_2, \frac{i}{m}, \frac{j}{m}\right)$ are involved. To find them by known ones, we will use the following relations

$$F\left(t, \frac{i}{m}H(\cdot-z_1), \frac{j}{m}H(\cdot-z_2)\right) = \iint_{\Delta_2} f\left(t, \xi_1, \xi_2, \frac{i}{m}H(\xi_1-z_1), \frac{j}{m}H(\xi_2-z_2)\right) d\xi_2 d\xi_1,$$

$$F\left(t, \frac{i}{m}H(\cdot-z_1), 0\right) = \int_0^{z_1}\int_0^{1-\xi_1} f(t, \xi_1, \xi_2, 0, 0) d\xi_2 d\xi_1 + \int_{z_1}^{1}\int_0^{1-\xi_1} f\left(t, \xi_1, \xi_2, \frac{i}{m}, 0\right) d\xi_2 d\xi_1,$$

$$F\left(t, 0, \frac{j}{m}H(\cdot-z_2)\right) = \int_0^{z_2}\int_0^{1-\xi_2} f(t, \xi_1, \xi_2, 0, 0) d\xi_1 d\xi_2 + \int_{z_2}^{1}\int_0^{1-\xi_2} f\left(t, \xi_1, \xi_2, 0, \frac{j}{m}\right) d\xi_1 d\xi_2,$$

from which we can easily obtain

$$f\left(t, z_1, z_2, \frac{i}{m}, \frac{j}{m}\right) = \frac{\partial^2}{\partial z_1 \partial z_2} F\left(t, \frac{i}{m}H(\cdot-z_1), \frac{j}{m}H(\cdot-z_2)\right) + f\left(t, z_1, z_2, \frac{i}{m}, 0\right) +$$

$$+ f\left(t, z_1, z_2, 0, \frac{j}{m}\right) - f(t, z_1, z_2, 0, 0),$$

$$\frac{\partial}{\partial z_1} F\left(t, \frac{i}{m}H(\cdot-z_1), 0\right) = \int_0^{1-z_1} f(t, z_1, \xi_2, 0, 0) d\xi_2 - \int_0^{1-z_1} f\left(t, z_1, \xi_2, \frac{i}{m}, 0\right) d\xi_2,$$

$$\frac{\partial}{\partial z_2} F\left(t, 0, \frac{j}{m} H(\cdot - z_2)\right) = \int_0^{1-z_2} f(t, \xi_1, z_2, 0, 0) d\xi_1 - \int_0^{1-z_2} f\left(t, \xi_1, z_2, 0, \frac{j}{m}\right) d\xi_1. \quad (15)$$

Let's substitute (15) into (14) and use the equalities

$$\sum_{i=0}^{m} \sum_{j=0}^{m-i} C_{i,j}(m) v_m^{(i,j)}(x, y, \alpha) g(i) = \sum_{i=0}^{m} C_m^i v_m^i(x, \alpha) g(i),$$

$$\sum_{i=0}^{m} \sum_{j=0}^{m-i} C_{i,j}(m) v_m^{(i,j)}(x, y, \alpha) g(j) = \sum_{j=0}^{m} \sum_{i=0}^{m-j} C_{i,j}(m) v_m^{(i,j)}(x, y, \alpha) g(j) = \sum_{j=0}^{m} C_m^j v_m^j(x, \alpha) g(j). \quad (16)$$

Then operator (14) can be written in the form

$$P_m^{[\alpha]}(F, x(\cdot), y(\cdot)) = \iint_{\Delta_2} \sum_{i=0}^{m} \sum_{j=0}^{m-i} C_{i,j}(m) v_m^{(i,j)}(x(z_1), y(z_2), \alpha) \frac{\partial^2}{\partial z_1 \partial z_2} F\left(t, \frac{i}{m} H(\cdot - z_1), \frac{j}{m} H(\cdot - z_2)\right) dz_1 dz_2 -$$
$$- \int_0^1 \sum_{i=0}^{m} C_m^i v_m^i(x(z_1), \alpha) \frac{\partial}{\partial z_1} F\left(t, \frac{i}{m} H(\cdot - z_1), 0\right) dz_1 - \int_0^1 \sum_{j=0}^{m} C_m^j v_m^j(y(z_2), \alpha) \frac{\partial}{\partial z_2} F\left(t, 0, \frac{j}{m} H(\cdot - z_2)\right) dz_2 +$$
$$+ F(t, 0, 0), \quad (17)$$

Where

$$v_m^{(i)}(x(z_1), \alpha) = \frac{\prod_{k_1=0}^{i-1}(x(z_1) + k_1 \alpha) \prod_{k_3=0}^{m-i-1}(1 - x(z_1) + k_3 \alpha)}{\prod_{k_4=0}^{m-1}(1 + k_4 \alpha)},$$

$$v_m^{(j)}(y(z_2), \alpha) = \frac{\prod_{k_1=0}^{j-1}(y(z_2) + k_1 \alpha) \prod_{k_3=0}^{m-j-1}(1 - y(z_2) + k_3 \alpha)}{\prod_{k_4=0}^{m-1}(1 + k_4 \alpha)}.$$

Formula (17) has a constructive character and it can be used in the practical calculations, but while theoretical investigating we will use formula (14) instead.

Let's consider, that

$$f(t, z_1, z_2, x, y) \in C([0,1]^5). \quad (18)$$

Then we will need

**Theorem 4** [3]. *If $f(x, y) \in C(\Delta)$ and $0 \leq \alpha = \alpha(m) \to 0$, when $m \to \infty$, then the sequence of operators $\{p_m^{[\alpha]}(f; x, y)\}$ uniformly converges to $f(x, y)$ on $\Delta$.*

Proving this theorem we used the next properties of operator polynomial $p_m^{[\alpha]}(f; x, y)$:

$$p_m^{[\alpha]}(1; x, y) = 1, \quad p_m^{[\alpha]}(t; x, y) = x, \quad p_m^{[\alpha]}(\tau; x, y) = y, \quad p_m^{[\alpha]}(t^2; x, y) = \frac{1}{1+\alpha}\left[\frac{x(1-x)}{m} + x(x+\alpha)\right],$$

$$p_m^{[\alpha]}(t\tau; x, y) = (1 - \frac{1}{m})\frac{xy}{1+\alpha}, \quad p_m^{[\alpha]}(\tau^2; x, y) = \frac{1}{1+\alpha}\left[\frac{y(1-y)}{m} + y(y+\alpha)\right]. \quad (19)$$

**Theorem 5.** *Let condition (18) holds and operator (13) is considered on the compact and* $0 \le \alpha = \alpha(m) \to 0$ *if* $m \to \infty$, *then the sequence of operators* $\{P_m^\alpha(F, x(\cdot), y(\cdot))\}$ *uniformly converges to* $F(t, x(\cdot), y(\cdot))$ *relatively to* $\{x(z_1), y(z_2)\} \in \Phi$, *where*
$\Phi = \{x(z_1) \in C[0,1]; y(z_2) \in C[0,1]; 0 \le x(z_1) + y(z_2) \le 1\}$.

*Proof.*

For every fixed $t, z_1, z_2 \in [0,1]$ operator polynomial $p_n^{[\alpha]}(f(t, z_1, z_2, \cdot, \cdot); x(z_1), y(z_2))$ from (14), in conformity with theorem 4 converges to $f(t, z_1, z_2, x(z_1), y(z_2))$ uniformly anywhere on the compact if $0 \le \alpha = \alpha(m) \to 0$ as $m \to \infty$. Then, from the evident generalization of theorem 1, p.506 from [8] (boundary transition in sub integral expression. Chapter. Integrals, that depend on the parameter), to our case the statement of theorem 1 follows.

### 6. ESTIMATE OF APPROXIMATION ORDER IN THE CASE OF TWO VARIABLES

For estimate of approximation order of operator $F(t, x(\cdot), y(\cdot))$ by operator polynomial (14) we use the modules of continuity, defined like that

$$\omega(\varphi, \delta) = \omega(\delta) = \max_{t, z_1, z_2 \in [0,1]} \sup_{|x'-x''|+|y'-y''| \le \delta} |\varphi(t, z_1, z_2, x'', y'') - \varphi(t, z_1, z_2, x', y')|,$$

where $\delta$ is positive number.

**Theorem 6.** *Let condition (18) hold and operator (13) is considered on the compact and* $\alpha \ge 0$, *then*

$$\left|F(t, x(\cdot), y(\cdot)) - P_n^{[\alpha]}(F, x(\cdot), y(\cdot))\right| \le 2\omega\left(\sqrt{\frac{1+\alpha m}{m + \alpha m}}\right).$$

*Proof.* So long as on the compact $\Phi = \{x(z_1) \in C[0,1]; y(z_2) \in C[0,1]; 0 \le x(z_1) + y(z_2) \le 1\}$ we have $C_{i,j}(m) v_m^{(i,j)}(x(z_1), y(z_2), \alpha) > 0$ and $P_m^{[\alpha]}(1, x(\cdot), y(\cdot)) = 1$, so, one can write

$$\left|F(t, x(\cdot), y(\cdot)) - P_n^{[\alpha]}(F, x(\cdot), y(\cdot))\right| \le$$

$$\le \iint_{\Delta_2} \sum_{i=0}^{m} \sum_{j=0}^{m-i} C_{i,j}(m) v_m^{(i,j)}(x(z_1), y(z_2), \alpha) \left|f(t, z_1, z_2, x(z_1), y(z_2)) - f\left(t, z_1, z_2, \frac{i}{m}, \frac{j}{m}\right)\right| dz_1 dz_2 \le$$

$$\le \iint_{\Delta_2} \sum_{i=0}^{m} \sum_{j=0}^{m-i} C_{i,j}(m) v_m^{(i,j)}(x(z_1), y(z_2), \alpha) \omega\left(f(t, z_1, z_2, \cdot, \cdot); \left|x(z_1) - \frac{i}{m}\right| + \left|\frac{j}{m} - y(z_2)\right|\right) dz_1 dz_2 \le$$

$$\le \iint_{\Delta_2} \sum_{i=0}^{m} \sum_{j=0}^{m-i} C_{i,j}(m) v_m^{(i,j)}(x(z_1), y(z_2), \alpha) \left(1 + \frac{1}{\delta}\left(\left|x(z_1) - \frac{i}{m}\right| + \left|\frac{j}{m} - y(z_2)\right|\right)\right) \omega(f(t, z_1, z_2, \cdot, \cdot); \delta) dz_1 dz_2 \le$$

$$\le \iint_{\Delta_2} \sum_{i=0}^{m} \sum_{j=0}^{m-i} C_{i,j}(m) v_m^{(i,j)}(x(z_1), y(z_2), \alpha) \left(1 + \frac{1}{\delta}\left(\left|x(z_1) - \frac{i}{m}\right| + \left|\frac{j}{m} - y(z_2)\right|\right)\right) dz_2 dz_1 \omega(\delta). \quad (20)$$

Here we have used the following properties of modules of continuity

$$|\varphi(x'', y'') - \varphi(x', y')| \le \omega(|x'' - x'| + |y'' - y'|); \quad \omega(\lambda \delta) \le (1 + \lambda)\omega(\delta) \text{ at }.$$

Using the identities (16), from (20) we obtain

$$\left|F(t, x(\cdot), y(\cdot)) - P_n^{[\alpha]}(F, x(\cdot), y(\cdot))\right| \le$$

$$\leq \omega(\delta) + \frac{1}{\delta}\int_0^1 \sum_{i=0}^m C_m^i v_m^{(i)}(x(z_1))\left|x(z_1) - \frac{i}{m}\right|dz_1 \omega(\delta) + \frac{1}{\delta}\int_0^1 \sum_{j=0}^m C_m^j v_m^{(j)}(y(z_2))\left|y(z_2) - \frac{j}{m}\right|dz_1 \omega(\delta).$$

Then, we act in the same way as while proving theorem 2. We have

$$\left|F(t, x(\cdot), y(\cdot)) - P_n^{[\alpha]}(F, x(\cdot), y(\cdot))\right| \leq \left(1 + \frac{1}{\delta}\sqrt{\frac{1+\alpha m}{m+\alpha m}}\right)\omega(\delta).$$

It is logically to choose $\delta = \sqrt{\dfrac{1+\alpha m}{m+\alpha m}}$, that leads to the completion of the theorem proof.

In case $\alpha = 0$ for Bernstein operator we will get

$$\left|F(t, x(\cdot), y(\cdot)) - B_m(F, x(\cdot), y(\cdot))\right| \leq 2\omega\left(\frac{1}{\sqrt{m}}\right).$$

## 7. ASYMPTOTIC ESTIMATE OF APPROXIMATION ERROR IN CASE OF TWO VARIABLES

Now we determine asymptotic estimate for the error

$$R_m^{[\alpha]}(F, x(\cdot), y(\cdot)) = F(t, x(\cdot), y(\cdot)) - P_m^{[\alpha]}(F, x(\cdot), y(\cdot)).$$

The following theorem is valid.

**Theorem 7.** *Let Urison operator (13) be considered on the compact and satisfy condition (18). For every fixed $t, z_1, z_2 \in [0,1]$ the second total differential $d^2 f(t, z_1, z_2, x, y)$ exists in the point $\{x(z_1), y(z_2)\} \in \Phi$, where $\Phi = \{x(z_1) \in C[0,1]; y(z_2) \in C[0,1]; 0 \leq x(z_1) + y(z_2) \leq 1\}$. Then, we have an asymptotic formula for the error*

$$R_m^{[\alpha]}(F, x(\cdot), y(\cdot)) =$$
$$= \iint_{\Delta_2} \sum_{i=0}^m \sum_{j=0}^{m-i} C_m^i C_{m-i}^j v_m^{(i,j)}(x(z_1), y(z_2), \alpha)\left(f(t, z_1, z_2, x(z_1), y(z_2)) - f(t, z_1, z_2, \frac{i}{m}, \frac{j}{m})\right)dz_1 dz_2 =$$
$$= -\frac{1+\alpha m}{1+\alpha}\iint_{\Delta_2}\left[\frac{x(z_1)(1-x(z_1))}{2m}f_{xx}''(t, z_1, z_2, x(z_1), y(z_2)) + \frac{x(z_1)y(z_2)}{m}f_{xy}''(t, z_1, z_2, x(z_1), y(z_2)) +\right.$$
$$\left. + \frac{y(z_2)(1-y(z_2))}{2m}f_{yy}''(t, z_1, z_2, x(z_1), y(z_2))\right]dz_1 dz_2 - \iint_{\Delta_2}\varepsilon_m^{[\alpha]}(x(z_1), y(z_2))dz_1 dz_2, \quad (21)$$

*where $m\varepsilon_m^{[\alpha]}(x(z_1), y(z_2))$ tends to zero, when $m \to \infty$.*

***Proof.*** Due to [3] and taking into account (19) we'll have

$$R_m^{[\alpha]}(F, x(\cdot), y(\cdot)) =$$
$$= \iint_{\Delta_2} \sum_{i=0}^m \sum_{j=0}^{m-i} C_m^i C_{m-i}^j v_m^{(i,j)}(x(z_1), y(z_2), \alpha)\left(f(t, z_1, z_2, x(z_1), y(z_2)) - f\left(t, z_1, z_2, \frac{i}{m}, \frac{j}{m}\right)\right)dz_1 dz_2 =$$
$$= -\frac{1+\alpha m}{1+\alpha}\iint_{\Delta_2}\left[\frac{x(z_1)(1-x(z_1))}{2m}f_{xx}''(t, z_1, z_2, x(z_1), y(z_2)) + \frac{x(z_1)y(z_2)}{m}f_{xy}''(t, z_1, z_2, x(z_1), y(z_2)) +\right.$$
$$\left. + \frac{y(z_2)(1-y(z_2))}{2m}f_{yy}''(t, z_1, z_2, x(z_1), y(z_2))\right]dz_1 dz_2 - \iint_{\Delta_2}\varepsilon_m^{[\alpha]}(x(z_1), y(z_2))dz_1 dz_2,$$

where

$$\varepsilon_m^{[\alpha]}(x(z_1), y(z_2)) = \frac{1}{2} \sum_{i=0}^{m} \sum_{j=0}^{m-i} C_m^i C_{m-i}^j v_m^{(i,j)}(x(z_1), y(z_2), \alpha) \times$$

$$\times \left\{ \left( x(z_1) - \frac{i}{m} \right)^2 \left[ \frac{\partial^2}{\partial x^2} f\left(t, z_1, z_2, x(z_1) + \theta\left(\frac{i}{m} - x(z_1)\right), y(z_2) + \theta\left(\frac{j}{m} - y(z_2)\right)\right) - \right.\right.$$

$$\left. - \frac{\partial^2}{\partial x^2} f(t, z_1, z_2, x(z_1), y(z_2)) \right] +$$

$$+ 2\left(\frac{i}{m} - x(z_1)\right)\left(\frac{j}{m} - y(z_2)\right) \left[ \frac{\partial^2}{\partial x \partial y} f\left(t, z_1, z_2, x(z_1) + \theta\left(\frac{i}{m} - x(z_1)\right), y(z_2) + \theta\left(\frac{j}{m} - y(z_2)\right)\right) - \right.$$

$$\left. - \frac{\partial^2}{\partial x \partial y} f(t, z_1, z_2, x(z_1), y(z_2)) \right] + \left(\frac{j}{m} - y(z_2)\right)^2 \times$$

$$\times \left[ \frac{\partial^2}{\partial y^2} f\left(t, z_1, z_2, x(z_1) + \theta\left(\frac{i}{m} - x(z_1)\right), y(z_2) + \theta\left(\frac{j}{m} - y(z_2)\right)\right) - \frac{\partial^2}{\partial y^2} f(t, z_1, z_2, x(z_1), y(z_2)) \right] \right\},$$

$$\theta \in (0,1).$$

From this relation we obtain the estimate

$$\left|\varepsilon_m^{[\alpha]}(x(z_1), y(z_2))\right| \leq \frac{1}{2} \sum_{i=0}^{m} \sum_{j=0}^{m-i} C_m^i C_{m-i}^j v_m^{(i,j)}(x(z_1), y(z_2), \alpha) \times$$

$$\times \left\{ \left( x(z_1) - \frac{i}{m} \right)^2 \left[ \omega\left( \frac{\partial^2}{\partial x^2} f(t, z_1, z_2, \cdot, \cdot); \left|x(z_1) - \frac{i}{m}\right| + \left|y(z_2) - \frac{j}{m}\right| \right) \right] + \right.$$

$$+ 2\left(\frac{i}{m} - x(z_1)\right)\left(\frac{j}{m} - y(z_2)\right)\left[ \omega\left( \frac{\partial^2}{\partial x \partial y} f(t, z_1, z_2, \cdot, \cdot); \left|x(z_1) - \frac{i}{m}\right| + \left|y(z_2) - \frac{j}{m}\right| \right) \right] +$$

$$\left. + \left(\frac{j}{m} - y(z_2)\right)^2 \left[ \omega\left( \frac{\partial^2}{\partial y^2} f(t, z_1, z_2, \cdot, \cdot); \left|x(z_1) - \frac{i}{m}\right| + \left|y(z_2) - \frac{j}{m}\right| \right) \right] \right\} \leq$$

$$\leq \frac{1}{2} \sum_{i=0}^{m} \sum_{j=0}^{m-i} C_m^i C_{m-i}^j v_m^{(i,j)}(x(z_1), y(z_2), \alpha) \left[ 1 + \frac{1}{\delta}\left( \left|x(z_1) - \frac{i}{m}\right| + \left|y(z_2) - \frac{j}{m}\right| \right) \right] \times$$

$$\times \left( \left|x(z_1) - \frac{i}{m}\right| + \left|y(z_2) - \frac{j}{m}\right| \right)^2 \omega(D^2 f; \delta),$$

where

$$\omega(D^2 f; \delta) = \max_{i=0,1,2} \max_{t, z_1, z_2 \in [0,1]} \omega(\frac{\partial^2}{\partial x^i \partial y^{2-i}} f(t, z_1, z_2, x, y); \delta).$$

Let's consider the limit of the last inequality, when $\alpha \to 0$, taking into account that functions $v_n^{(k)}(x(z_1), y(z_2), \alpha)$ are continuous with respect to $\alpha$,. Then, we have

$$\left|\varepsilon_m^{[0]}(x(z_1), y(z_2))\right| \leq \frac{1}{2} \sum_{i=0}^{m} \sum_{j=0}^{m-i} C_m^i C_{m-i}^j x(z_1)^i y(z_2)^j (1 - x(z_1) - y(z_2))^{m-i-j} \times$$

$$\times \left[ 1 + \frac{1}{\delta}\left(\left|x(z_1) - \frac{i}{m}\right| + \left|y(z_2) - \frac{j}{m}\right|\right)\right] \left(\left|x(z_1) - \frac{i}{m}\right| + \left|y(z_2) - \frac{j}{m}\right|\right)^2 \omega(D^2 f; \delta) \le$$

$$\le \frac{1}{2}\left\{\sum_{i=0}^{m}\sum_{j=0}^{m-i} C_m^i C_{m-i}^j x(z_1)^i y(z_2)^j (1 - x(z_1) - y(z_2))^{m-i-j} \left(\left|x(z_1) - \frac{i}{m}\right| + \left|y(z_2) - \frac{j}{m}\right|\right)^4\right\}^{1/2} \times$$

$$\times \left\{\sum_{i=0}^{m}\sum_{j=0}^{m-i} C_m^i C_{m-i}^j x(z_1)^i y(z_2)^j (1 - x(z_1) - y(z_2))^{m-i-j} \left(1 + \frac{1}{\delta}\left(\left|x(z_1) - \frac{i}{m}\right| + \left|y(z_2) - \frac{j}{m}\right|\right)\right)^2\right\}^{1/2} \omega(D^2 f; \delta).$$

Then, for the first root from the first right side of the last inequality we will use the inequality $(U_1 + U_2)^2 \le 2(U_1^2 + U_2^2)$, twice, and for the second one we will use the inequality $\left(\sum_{i=1}^{m} U_i\right)^2 \le m\sum_{i=1}^{m} U_i^2$. Then, we obtain

$$\left|\varepsilon_m^{[0]}(x(z_1), y(z_2))\right| \le \sqrt{3}\left\{\sum_{i=0}^{m}\sum_{j=0}^{m-i} C_m^i C_{m-i}^j x(z_1)^i y(z_2)^j (1 - x(z_1) - y(z_2))^{m-i-j}\right.$$

$$\left.\left(\left|x(z_1) - \frac{i}{m}\right|^4 + \left|y(z_2) - \frac{j}{m}\right|^4\right)\right\}^{1/2} \times$$

$$\times \left\{\sum_{i=0}^{m}\sum_{j=0}^{m-i} C_m^i C_{m-i}^j x(z_1)^i y(z_2)^j (1 - x(z_1) - y(z_2))^{m-i-j} \left(1 + \frac{1}{\delta^2}\left(\left|x(z_1) - \frac{i}{m}\right|^2 + \left|y(z_2) - \frac{j}{m}\right|^2\right)\right)\right\}^{1/2} \times$$

$$\times \omega(D^2 f; \delta) \le \sqrt{3}\left(\frac{1}{m^3}\left(\frac{3m}{8} - \frac{1}{4}\right)\right)^{1/2}\left(1 + \frac{1}{2\delta^2 m}\right)^{1/2} \omega(D^2 f; \delta) \le \frac{3\sqrt{3}}{4m}\omega\left(D^2 f; \frac{1}{\sqrt{m}}\right).$$

To obtain this estimate the identities (11), (19) and value $\delta = \frac{1}{\sqrt{m}}$ have been used. So, when $m$, is big enough the following inequality will hold

$$\left|\varepsilon_m^{[0]}(x(z_1), y(z_2))\right| \le \frac{3\sqrt{3}}{4m}\omega\left(D^2 f; \frac{1}{\sqrt{m}}\right).$$

Let's denote $m\varepsilon_m^{[\alpha]}(x(z_1), y(z_2)) = \frac{3\sqrt{3}}{4m}\omega\left(D^2 f; \frac{1}{\sqrt{m}}\right)$. So, the theorem has been proved.

In particular case, when $\alpha = 0$, for Bernstein polynomial we'll obtain an asymptotic formula

$$F(t, x(\cdot), y(\cdot)) = B_n(F, x(\cdot), y(\cdot)) -$$

$$- \iint_{\Delta_2}\left[\frac{x(z_1)(1 - x(z_1))}{2m}f_{xx}''(t, z_1, z_2, x(z_1), y(z_2)) + \frac{x(z_1)y(z_2)}{m}f_{xy}''(t, z_1, z_2, x(z_1), y(z_2)) +\right.$$

$$\left. + \frac{y(z_2)(1 - y(z_2))}{2m}f_{yy}''(t, z_1, z_2, x(z_1), y(z_2))\right]dz_1 dz_2 - \iint_{\Delta_2}\varepsilon_m^{[0]}(x(z_1), y(z_2))dz_1 dz_2.$$

Let $\alpha = -\dfrac{1}{m}$. Whereas in this case polynomial $P_m^{\left[-\frac{1}{m}\right]}(F, x(\cdot), y(\cdot))$ keeps the polynomial of two variables of $m$ - degree, then substituting Taylor series into (21) for $f(t, z_1, z_2, \dfrac{i}{m}, \dfrac{j}{m})$

$$f(t, z_1, z_2, \frac{i}{m}, \frac{j}{m}) = \sum_{p=0}^{\infty} \frac{(-1)^p}{p!} \left( (x(z_1) - \frac{i}{m}) \frac{\partial}{\partial x} + (y(z_2) - \frac{j}{m}) \frac{\partial}{\partial y} \right)^p f(t, z_1, z_2, x(z_1), y(z_2)),$$

we obtain

$$P_m^{\left[-\frac{1}{m}\right]}(F, x(\cdot), y(\cdot)) = F(t, x(\cdot), y(\cdot)) + \sum_{n=m}^{\infty} \frac{(-1)^{n+1}}{(n+1)!} \iint_{\Delta_2} \sum_{i=0}^{m} \sum_{j=0}^{m-i} C_m^i C_{m-i}^j v_m^{(i,j)}\left( x(z_1), y(z_2), -\frac{1}{m} \right) \times$$

$$\times \left( (x(z_1) - \frac{i}{m}) \frac{\partial}{\partial x} + (y(z_2) - \frac{j}{m}) \frac{\partial}{\partial y} \right)^{n+1} f(t, z_1, z_2, x(z_1), y(z_2)) dz_2 dz_1. \qquad (22)$$

It is impossible to make reasoning similar to one-dimensional case, because here the polynomials $v_m^{(i,j)}\left( x(z_1), y(z_2), -\dfrac{1}{m} \right)$, are not nonnegative . Therefore, let's use the following designation

$$a_{k,n+1-k}^{\left[-\frac{1}{m}\right]}\left( x(z_1), y(z_2), -\frac{1}{m} \right) = \sum_{i=0}^{m} \sum_{j=0}^{m-i} C_m^i C_{m-i}^j v_m^{(i,j)}(x(z_1), y(z_2), -\frac{1}{m})[x(z_1) - \frac{i}{m}]^k [y(z_2) - \frac{j}{m}]^{n+1-k}$$

to estimate of behavior of

$$-R_n^{\left[-\frac{1}{m}\right]}(F, x(\cdot), y(\cdot)) = P_m^{\left[-\frac{1}{m}\right]}(F, x(\cdot), y(\cdot)) - F(t, x(\cdot), y(\cdot)).$$

Then, formula (22) takes the form

$$P_m^{\left[-\frac{1}{m}\right]}(F, x(\cdot), y(\cdot)) = F(t, x(\cdot), y(\cdot)) +$$
$$+ \sum_{n=m}^{\infty} \frac{(-1)^{n+1}}{(n+1)!} \iint_{\Delta_2} \sum_{k=0}^{n+1} C_{n+1}^k a_{k,n+1-k}^{\left[-\frac{1}{m}\right]}\left( x(z_1), y(z_2), -\frac{1}{m} \right) \frac{\partial^{n+1}}{\partial x^k \partial y^{n+1-k}} f(t, z_1, z_2, x(z_1), y(z_2)) dz_2 dz_1. \qquad (23)$$

One can verify that following correlations occurs

$$A_{n+1}^{\left[-\frac{1}{m}\right]} = \max_{\substack{x,y \geq 0 \\ x+y \leq 1}} \left| a_{k,n+1-k}^{\left[-\frac{1}{m}\right]}\left( x, y, -\frac{1}{m} \right) \right| =$$

$$= \max_{0 \leq x \leq 1} \left| a_{k,n+1-k}^{\left[-\frac{1}{m}\right]}\left( x, 1-x, -\frac{1}{m} \right) \right| = \max_{0 \leq x \leq 1} \left| a_{n+1,0}^{\left[-\frac{1}{m}\right]}\left( x, 1-x, -\frac{1}{m} \right) \right|, \quad k = 0,1,...,n+1, \qquad (24)$$

$$\left| a_{k,n+1-k}^{\left[-\frac{1}{m}\right]}\left( x, 1-x, -\frac{1}{m} \right) \right| = \frac{m^{m-1}}{(m-1)!} \left| x(x - \frac{1}{m})...(x - \frac{m}{m}) \sum_{i=0}^{m} (-1)^i C_m^i (x - \frac{i}{m})^{n+1} \right|, \quad k = 0,1,...,n+1.$$

Taking into account (24) from (23) we obtain the inequality

$$\left| R_n^{\left[-\frac{1}{m}\right]}(F, x(\cdot), y(\cdot)) \right| \le$$

$$\le \sum_{n=m}^{\infty} \frac{1}{(n+1)!} \iint_{\Delta_2} \left| \sum_{k=0}^{n+1} C_{n+1}^k a_{k,n+1-k}^{\left[-\frac{1}{m}\right]}\left(x(z_1), y(z_2), -\frac{1}{m}\right) \right| \left| \frac{\partial^{n+1}}{\partial x^k \partial y^{n+1-k}} f(t, z_1, z_2, x(z_1), y(z_2)) \right| dz_2 dz_1 \le$$

$$\le \sum_{n=m}^{\infty} \frac{2^{n+1}}{(n+1)!} A_{n+1}^{\left[-\frac{1}{m}\right]} D^{n+1} f, \qquad (25)$$

where

$$D^{n+1} f = \max_{\substack{0 \le t, z_1, z_2 \le 1 \\ z_1 + z_2 \le 1}} \max_{\substack{0 \le x(z_1), y(z_2) \le 1 \\ x(z_1) + y(z_2) \le 1}} \left| \frac{\partial^{n+1}}{\partial x^k \partial y^{n+1-k}} f(t, z_1, z_2, x(z_1), y(z_2)) \right|.$$

Let's consider the first (main) summand in a right-hand side of inequality (25). As,

$$\left| a_{k,m+1-k}^{\left[-\frac{1}{m}\right]}\left(x, 1-x, -\frac{1}{m}\right) \right| = \frac{m^{m-1}}{(m-1)!} \left| x(x-\frac{1}{m})\ldots(x-\frac{m}{m}) \sum_{i=0}^{m} (-1)^i C_m^i (x-\frac{i}{m})^{m+1} \right| =$$

$$= (m+1) \left| (x-\frac{1}{2}) \, x(x-\frac{1}{m})\ldots(x-\frac{m}{m}) \right|$$

then

$$A_{m+1}^{\left[-\frac{1}{m}\right]} = (m+1) \max_{x \in [0,1]} \left| (x-\frac{1}{2}) \prod_{i=0}^{m} (x-\frac{i}{m}) \right| = \frac{m+1}{m^{m+2}} \max_{t \in [0,m]} \left| (t-\frac{m}{2}) \prod_{i=0}^{m} (t-i) \right| = \frac{m+1}{m^{m+2}} \max_{t \in [0,m]} |\varphi(t)|. \qquad (26)$$

To estimate a right-hand side of relation (26) analogously to [9], p.95, we'll consider the function

$$\varphi(t) = (t-\frac{m}{2}) \prod_{i=0}^{m} (t-i) = \varphi(z+\frac{m}{2}) = z\left[ z^2 - \left(\frac{m}{2}\right)^2 \right]\left[ z^2 - \left(\frac{m-2}{2}\right)^2 \right],$$

that as the function of $z$ is even or odd with respect to evenness or oddness of $m$. The following relation is valid

$$\varphi(t+1) = v(t)\varphi(t), \quad v(t) = \frac{t+1}{t-m} \frac{t+1-\frac{m}{2}}{t-\frac{m}{2}},$$

where function $v(t)$ is negative, when $t$ changes from $0$ to $\frac{m}{2}-1$ and its module comes up to maximum at point

$$t_{max} = \frac{m-1-\sqrt{m+1}}{2},$$

that can be defined by the formula

$$\max_{0 \leq t \leq \frac{m}{2}-1} |v(t)| = |v(t_{\max})| = \left[\frac{-1+\sqrt{m+1}}{1+\sqrt{m+1}}\right]^2. \qquad (27)$$

From (27) we obtain, that $|v(t_{\max})|$ is an increasing function relatively to $m \in \{3,4,...\}$, its maximum is reached in the infinity and is equal to 1. Thus, the extreme values of $\varphi(t)$ will decrease in module up to the middle of the interval $[0,m]$, and then by symmetry (antisymmetry) will increase. Above-stated give grounds to conclude: there exists the point $x = \xi \in (0, \frac{1}{m})$, where the equality holds

$$A_{m+1}^{[-\frac{1}{m}]} = (m+1)\left|(\xi - \frac{1}{2})\prod_{i=0}^{m}(\xi - \frac{i}{m})\right|. \qquad (28)$$

Hence, the estimate follows

$$A_{m+1}^{[-\frac{1}{m}]} \leq \frac{1}{2}\frac{(m+1)!}{m^{m+1}}.$$

Then we estimate $A_{m+2}^{[-\frac{1}{m}]}$. So, we have

$$A_{m+2}^{[-\frac{1}{m}]} = \frac{m^m}{m!}\max_{0 \leq x \leq 1}\left|x(x-\frac{1}{m})...(x-\frac{m}{m})\sum_{i=0}^{m}(-1)^i C_m^i (x-\frac{i}{m})^{m+2}\right| =$$

$$= \frac{1}{m!m^{m+3}}\max_{0 \leq t \leq m}\left|\omega(t)\sum_{i=0}^{m}(-1)^i C_m^i (t-i)^{m+2}\right| = \max_{0 \leq t \leq m}\left|a_{m+2,0}^{[-\frac{1}{m}]}(\frac{t}{m},1-\frac{t}{m},-\frac{1}{m})\right|,$$

where $\omega(t) = t(t-1)(t-2)...(t-m)$. The following relation is valid

$$a_{m+2,0}^{[-\frac{1}{m}]}(\frac{t+1}{m},1-\frac{t+1}{m},-\frac{1}{m}) = v_{m+2}(t)a_{m+2,0}^{[-\frac{1}{m}]}(\frac{t}{m},1-\frac{t}{m},-\frac{1}{m}),$$

$$v_{m+2}(t) = \frac{t+1}{t-m}\frac{\sum_{i=0}^{m}(-1)^i C_m^i (t+1-i)^{m+2}}{\sum_{i=0}^{m}(-1)^i C_m^i (t-i)^{m+2}}. \qquad (29)$$

And the function $v_{m+2}(t)$ meets the conditions

$$v_{m+2}(t) < 0, \ t \in [0, \frac{m-1}{2}], \quad \max_{t \in [0, \frac{m-1}{2}]}|v_{m+2}(t)| = \left|v_{m+2}(\frac{m-1}{2})\right| = 1, \qquad (30)$$

that can be evident by means of direct verification.

Rely on foresaid reasoning one can say that the point $\xi \in (0, \frac{1}{m})$ exists and the following correlation is valid

$$A_{m+2}^{[-\frac{1}{m}]} = \frac{m^m}{m!}\left|\xi(\xi-\frac{1}{m})...(\xi-\frac{m}{m})\sum_{i=0}^{m}(-1)^i C_m^i (\xi-\frac{i}{m})^{m+2}\right| \le 2^m.$$

Let's show, that the same estimate will occur also for $A_{n+1}^{[-\frac{1}{m}]}$, $\forall n = m, m+1, m+2,...$ . We need to introduce the next function

$$f_{m,n}(t) = \sum_{i=0}^{m}(-1)^i C_m^i (t-i)^{n+1},$$

that is followed by a valid correlation

$$\frac{d}{dt}f_{m,n}(t) = (n+1)f_{m,n-1}(t). \tag{31}$$

Then we use a certain integral representation of difference of arbitrary order (see, for instance [9]), that has the consequence in a form

$$f_{m,n}(t) = \sum_{i=0}^{m}(-1)^i C_m^i (t-i)^{n+1} = \frac{m!(n+1)!}{(n-m+1)!}\int_0^1\int_0^{z_1}...\int_0^{z_{m-1}}(t-z_1-z_2-...-z_m)^{n-m+1}dz_m...dz_1 \tag{32}$$

Let $m+n$ be even number, then from (32) follows that the function $f_{m,n}(t)$ is increasing in the interval $[0, m]$. But since in this case

$$f_{m,n}(\frac{m}{2}) = 0, \tag{33}$$

then $f_{m,n}(0) < f_{m,n}(t) < 0$, $t \in (0, \frac{m}{2})$. This implies, that

$$0 < \frac{f_{m,n}(t+1)}{f_{m,n}(t)} < 1, \quad t \in [0, \frac{m}{2}-1]. \tag{34}$$

The validity of (33) follows from the equality

$$f_{m,n}(\frac{m}{2}) = \left(\frac{m}{2}\right)^{n+1} + (-1)^m\left(-\frac{m}{2}\right)^{n+1} + C_m^1\left(-\left(\frac{m}{2}-1\right)^{n+1} + (-1)^{m-1}\left(-\frac{m}{2}+1\right)^{n+1}\right) + ... = 0$$

Now let $m+n$ be odd number, then the function $f_{m,n}(t)$ is a decreasing function in the interval $(0, \frac{m}{2})$ and the inequalities $0 < f_{m,n}(\frac{m}{2}) < f_{m,n}(t) < f_{m,n}(0)$, $t \in (0, \frac{m}{2})$ is valid. They also lead to statement (34), as in previous case. Established inequalities (34) lead to the conclusion, that the function

$$v_{n+1}(t) = \frac{t+1}{t-m}\frac{f_{m,n}(t+1)}{f_{m,n}(t)} = \frac{t+1}{t-m}\frac{\sum_{i=0}^{m}(-1)^i C_m^i (t+1-i)^{n+1}}{\sum_{i=0}^{m}(-1)^i C_m^i (t-i)^{n+1}}$$

is negative and module of it is less than 1 in the interval $(0, \frac{m}{2}-1)$. Then we have

$$A_{n+1}^{[-\frac{1}{m}]} = (m+1)\max_{x\in[0,1]}\left|\prod_{i=0}^{m}(x-\frac{i}{m})\sum_{i=0}^{m}(-1)^i C_m^i (x-\frac{i}{m})^{n+1}\right| = \frac{m+1}{m^{m+2}}\max_{t\in[0,m]}\left|\prod_{i=0}^{m}(t-i)f_{m,n}(t)\right| = \frac{m+1}{m^{m+2}}\max_{t\in[0,m]}|\phi(t)|.$$

Hence, since, $\varphi(t+1) = \frac{t+1}{t-m}\frac{f_{m,n}(t+1)}{f_{m,n}(t)}\varphi(t)$, then taking into account (34) one can state, that maximum modulo values of function $\varphi(t)$ decreases in the interval $\left[0, \frac{m}{2}\right]$. So, there exists the point like $\xi \in (0, \frac{1}{m})$, that the following relation is valid

$$A_{n+1}^{[-\frac{1}{m}]} = \frac{m^m}{m!}\left|\xi(\xi-\frac{1}{m})...(\xi-\frac{m}{m})\sum_{i=0}^{m}(-1)^i C_m^i (\xi-\frac{i}{m})^{n+1}\right| \le 2^m, \forall n = m, m+1,... \qquad (35)$$

The inequality (35) together with (25) lead to the estimate

$$|R_n^{[-\frac{1}{m}]}(F, x(\cdot), y(\cdot))| \le 2^m \sum_{n=m}^{\infty} \frac{2^{n+1}}{(n+1)!} D^{n+1}f,$$

which is followed by

**Theorem 8.** *Let function* $f(t, z_1, z_2, x, y)$ *be such, that in the domain* $0 \le t \le 1$ $0 \le z_1, z_2, z_1+z_2 \le 1$, $0 \le x, y \le 1$, $x+y \le 1$ *there exist* $D^{n+1}f$, $\forall n = 1,2,...$ *such, the series*

$$\frac{D^{m+1}f}{m+1} + \frac{2D^{m+2}f}{(m+1)(m+2)} + \frac{2^2 D^{m+3}f}{(m+1)(m+2)(m+3)} + ...$$

*is convergent and its sum has upper estimate* $M$, *that doesn't depend on* $m$. *Then, the following inequalities are valid*

$$|R_n^{[-\frac{1}{m}]}(F, x(\cdot), y(\cdot))| \le \frac{2^{2m+1}}{m!}M \le \frac{\sqrt{2}}{\sqrt{\pi n}}\left(\frac{4e}{m}\right)^m M.$$

*Example 2.* Let's consider Urison operator

$$F(t, x(\cdot), y(\cdot)) = \iint_{\Delta_2} \sin(x(z_1))\cos(y(z_2))\,dz_1 dz_2, \quad \{x(z_1), y(z_2)\} \in \Phi$$

and construct operator polynomial (14) for it

$$P_m^{[\alpha]}(F, x(\cdot), y(\cdot)) = \iint_{\Delta_2} \sum_{0 \le i+j \le m} \frac{m!}{i!\,j!(m-i-j)!} v_m^{(i,j)}(x(z_1), y(z_2), \alpha)\sin(x(z_1))\cos(y(z_2))dz_1 dz_2,$$

where $v_m^{(i,j)}(x(z_1), y(z_2), \alpha) = \dfrac{\prod\limits_{k_1=0}^{i-1}(x(z_1)+k_1\alpha)\prod\limits_{k_2=0}^{j-1}(y(z_2)+k_2\alpha)\prod\limits_{k_3=0}^{m-i-j-1}(1-x(z_1)-y(z_2)+k_3\alpha)}{\prod\limits_{k_4=0}^{m-1}(1+k_4\alpha)}$, $\alpha \geq 0$.

We choose, for example, and $\alpha = -\dfrac{1}{m}$ i $x(z_1) = z_1/(1+z_1^2)$, $y(z_2) = z_2/(1+z_2)$. For calculation we'll use Maple. The results we'll write in table 2, where
$\Delta_1 = |F(t, x(\cdot), y(\cdot)) - P_m^{[1/m]}(F, x(\cdot), y(\cdot))|$, $\Delta_2 = |F(t, x(\cdot), y(\cdot)) - B_m(F, x(\cdot), y(\cdot))|$,
$\Delta_3 = |F(t, x(\cdot), y(\cdot)) - P_m^{[-1/m]}(F, x(\cdot), y(\cdot))|$.

Table 2

| n | $\Delta_1$ | $C_1 = m*\Delta_1$ | $\Delta_2$ | $C_2 = m*\Delta_2$ | $\Delta_3$ |
|---|---|---|---|---|---|
| 1 | 0.013248971 | 0.013248972 | 0.1324897e-1 | 0.013248972 | 0.013248972 |
| 2 | 0.87424439e-2 | 0.017484888 | 0.6489178e-2 | 0.012978357 | 0.0002706 |
| 4 | 0.52128810e-2 | 0.020851524 | 0.3232767e-2 | 0.012931069 | 0.236968e-6 |
| 8 | 0.28857948e-2 | 0.023086359 | 0.16160302e-2 | 0.012928242 | 0.17266e-13 |
| 16 | 0.15249613e-2 | 0.024399381 | 0.8080102e-3 | 0.012928163 | 0.53346e-26 |

We see, that the inequalities

$$\Delta_1 \leq \frac{0.03}{m}, \quad \Delta_2 \leq \frac{0.015}{m}, \quad \Delta_3 \leq \frac{2^{2m}}{m!}(e^2-1), \quad (M < \frac{1}{2}(e^2-1))$$

are valid.

***Example 3.*** Let's consider Urison operator

$$F(t, x(\cdot), y(\cdot)) = \iint_{\Delta_2}(1 - 0{,}5z_2 x - 0{,}5z_1 y)^4 \, dz_1 dz_2, \quad \{x(z_1), y(z_2)\} \in \Phi$$

We construct the operator polynomial (14) for it

$$P_m^{[\alpha]}(F, x(\cdot), y(\cdot)) = \iint_{\Delta_2} \sum_{0 \leq i+j \leq m} \frac{m!}{i!\,j!\,(m-i-j)!} v_m^{(i,j)}(x(z_1), y(z_2), \alpha)\left(1 - 0{,}5z_2\frac{i}{m} - 0{,}5z_1\frac{1}{m}\right)^4 dz_1 dz_2,$$

where $v_m^{(i,j)}(x(z_1), y(z_2), \alpha) =$

$= \dfrac{\prod\limits_{k_1=0}^{i-1}(x(z_1)+k_1\alpha)\prod\limits_{k_2=0}^{j-1}(y(z_2)+k_2\alpha)\prod\limits_{k_3=0}^{m-i-j-1}(1-x(z_1)-y(z_2)+k_3\alpha)}{\prod\limits_{k_4=0}^{m-1}(1+k_4\alpha)}$, $\alpha \geq 0$.

Let's choose, for example, $\alpha = \dfrac{1}{m}$, $\alpha = 0$ and $\alpha = -\dfrac{1}{m}$ i $x(z_1) = z_1/(1+z_1^2)$, $y(z_2) = z_2/(1+z_2)$. For calculation we use Maple. The results we write in Table 3, where
$\Delta_1 = |F(t, x(\cdot), y(\cdot)) - P_m^{[1/m]}(F, x(\cdot), y(\cdot))|$, $\Delta_2 = |F(t, x(\cdot), y(\cdot)) - B_m(F, x(\cdot), y(\cdot))|$,
$\Delta_3 = |F(t, x(\cdot), y(\cdot)) - P_m^{[-1/m]}(F, x(\cdot), y(\cdot))|$.

Table 3

| n  | $\Delta_1$ | $C_1 = m*\Delta_1$ | $\Delta_2$ | $C_2 = m*\Delta_2$ | $\Delta_3$ |
|----|----------|-------------------|----------|-------------------|----------|
| 1  | 0.016387 | 0.016387          | 0.016387 | 0.016387          | 0.016387 |
| 2  | 0.011145 | 0.022290          | 0.008524 | 0.017048          | 0.00066  |
| 4  | 0.006819 | 0.027274          | 0.004347 | 0.017386          | 0        |
| 8  | 0.003848 | 0.030780          | 0.002195 | 0.017557          | 0        |
| 16 | 0.002058 | 0.032934          | 0.001103 | 0.0176430         | 0        |
| 32 | 0.001067 | 0.034143          | 0.000553 | 0.017686          | 0        |

We see, that inequalities

$$\Delta_1 \leq \frac{0.036}{m}, \quad \Delta_2 \leq \frac{0.018}{m},$$

are valid and beginning from $n = 4$ the approximation error $\Delta_3$ is equal to zero.

## 8. CONCLUSIONS

In this work we investigate a linear polynomial operator that approximates Urison operator, when integration domain is a "regular triangle". The basis of this polynomial are the polynomials, suggested by D.D. Stancu and from which Bernstein polynomials are obtained. Theoretical and computation investigations showed, that under the high kernel smoothness of Urison operator, the best accuracy is obtained by suggested in this work operator polynomial with the parameter $\alpha = -\frac{1}{m}$ in comparison to those, in which $\alpha \geq 0$. Here $m$ is a degree of operator polynomial. Note, that this polynomial has a wonderful characteristic: it interpolates Urison operator on the continuous set of knots.